\documentclass[11pt]{amsart}
\usepackage{amsfonts,amscd,amsmath,amssymb,amsthm,mathrsfs}
\usepackage{array,latexsym,color,float,enumerate}
\usepackage{graphicx,epsfig}
\usepackage[all]{xy}
\usepackage{pdfsync}
\usepackage{marginnote}
\usepackage[hypertex,linkbordercolor={0 0 1}]{hyperref}

\theoremstyle{plain}
\newtheorem{theorem}{Theorem}[section]
\newtheorem{lemma}[theorem]{Lemma}
\newtheorem{proposition}[theorem]{Proposition}
\newtheorem{corrolary}[theorem]{Corollary}
\theoremstyle{definition}
\newtheorem{definition}[theorem]{Definition}
\newtheorem{example}[theorem]{Example}
\newtheorem*{conjecture}{Conjecture}

\newtheorem{remark}[theorem]{Remark}
\newtheorem{construction}[theorem]{Construction}
\theoremstyle{remark}

\newcommand{\bthm}{\begin{theorem}}
\newcommand{\bprop}{\begin{proposition}}
\newcommand{\blem}{\begin{lemma}}
\newcommand{\bcor}{\begin{corrolary}}
\newcommand{\brem}{\begin{remark}}
\newcommand{\bdfn}{\begin{definition}}
\newcommand{\bex}{\begin{example}}
\newcommand{\bcon}{\begin{construction}}
\newcommand{\econ}{\end{construction}}
\newcommand{\eex}{\end{example}}
\newcommand{\ethm}{\end{theorem}}
\newcommand{\eprop}{\end{proposition}}
\newcommand{\elem}{\end{lemma}}
\newcommand{\ecor}{\end{corrolary}}
\newcommand{\erem}{\end{remark}}
\newcommand{\edfn}{\end{definition}}
\newcommand{\benum}{\begin{enumerate}}
\newcommand{\eenum}{\end{enumerate}}

\newcommand{\ov}{\overline}

\renewcommand{\frak}[1]{\mathfrak{#1}}

\def\8{\infty}
\def\.{\cdot}
\def\PP{\mathbb{P}}
\def\F{\mathbb{F}}
\def\C{\mathbb{C}}
\def\Z{\mathbb{Z}}

\def\Q{\mathbb{Q}}

\def\E{\widehat{E}}
\def\ovk{\overline\kappa}
\def\kp{\kappa}

\def\:{\colon}
\def\map{\dashrightarrow}

\def\ssk{\smallskip}

\newcommand{\noin}{\noindent}
\def\Bk{\operatorname{Bk}}

\def\Aut{\operatorname{Aut}}

\def\Coker{\operatorname{Coker}}

\def\im{\operatorname{Im}}

\def\Pic{\operatorname{Pic}}
\def\dim{\operatorname{dim}}

\begin{document}
\title[Recent progress in the geometry of $\Q$-acyclic surfaces]{Recent progress in the geometry of $\Q$-acyclic surfaces}

\author[Karol Palka]{Karol Palka}\address{Karol Palka: Department of Mathematics, McGill University,
Montreal, QC, Canada}\email{palka@impan.pl}

\thanks{The author was supported by Polish Grant MNiSzW}


\begin{abstract} \noindent We give a survey of results on the geometry of complex algebraic  $\Q$-acyclic surfaces including some recent results. \end{abstract}

\maketitle

\tableofcontents

\section{Introduction}

In this article we want to take the reader to the rich and yet not fully explored world of plane-like complex algebraic surfaces. We hope our survey will give a general reader a taste of methods used and will serve as an update on recent results for experts. All varieties considered are complex algebraic.

The story begins with a surprising discovery of Ramanujam \cite{Ramanujam} of a contractible affine surface non-isomorphic to $\C^2$, one of the many nontrivial \emph{smooth homotopy planes}, i.e. smooth contractible surfaces. Ramanujam discovered a first important characterization of $\C^2$ by showing that it is the unique smooth homotopy plane which is simply connected at infinity. Since then affine algebraic geometers began the study of smooth (and more generally normal) varieties with the same Betti numbers as $\C^2$, so-called \emph{$\Q$-homology planes}. One of the motivations was the search for more practical characterizations of the complex plane (the computation of the fundamental group at infinity of an affine variety is usually very difficult). In the face of the topological simplicity the intriguing algebraic side of $\Q$-acyclic surfaces is more clearly visible. Because of their homological similarity to the plane, smooth $\Q$-homology planes, and especially smooth homotopy planes, play today an important role as a source of examples or counterexamples when studying working hypotheses as well as more challenging conjectures. They accompany us when studying exotic structures on $\C^n$'s (see \cite{Zaid_exotic_lectures}), the Cancellation Conjecture\footnote{The $m$-dimensional Cancellation Conjecture: If $X\times \C^n\cong \C^{n+m}$ then $X\cong \C^m$.} (which was the motivation for Ramanujam and has been proved in dimension two by Fujita and Miyanishi \cite{Fuj-Zar}), the Jacobian Conjecture\footnote{The Jacobian Conjecture: A polynomial map $f\:\C^n\to\C^n$ with nowhere vanishing Jacobian determinant $|J(f)|$ has a polynomial inverse.} (see \cite[5.2]{Miyan_recent_developments}) and others. Today, after almost forty years, mainly due to the tools of the theory of open surfaces, the scheme of the classification, and in most cases the classification itself, are settled. Recently the author had the pleasure of adding his part to the story.

Although we try to avoid notions which are not well-known to any well-versed algebraic geometer, we need to refer the reader for the basic notions of the theory of open surfaces to \cite{Miyan-OpenSurf}. Since the section 4.3 loc. cit. is a review of smooth $\Q$-homology planes we concentrate mainly on the singular case, stating the updated results so to include the smooth case if possible (see also \cite{Miyan_recent_developments}).

\ssk\textsl{\textsf{Acknowledgements.}} The author would like to thank professors P. Russell and S. Lu for the invitation to Montreal. He also thanks prof. P. Russell for stimulating discussions.

\section{Preliminaries}

A \emph{normal pair} $(\ov X,D)$ consists of a complete normal surface $\ov X$ and a divisor $D$ with simple normal crossings (snc) contained in the smooth locus of $\ov X$. A normal pair $(\ov X,D)$ is said to be smooth if $\ov X$ is smooth. If $X$ is a normal surface then an embedding $\iota\:X\to \ov X$, where $(\ov X,\ov X\setminus X)$ is a normal pair, is called a \emph{normal completion} of $X$. A normal completion is a \emph{smooth completion} if $\ov X$ is smooth. Two normal completions $\iota_j\:X\to \ov X_j$, $j=1,2$ are isomorphic if there exists a morphism $f:\ov X_1\to\ov X_2$ for which $f\circ\iota_1=\iota_2$. A surface with isolated singularities is called \emph{logarithmic} if all its singular points are of analytical type $\C^2/G$ for some finite subgroup $G<GL(2,\C)$. By an \emph{$n$-curve} we mean a smooth rational curve with self-intersection $n$.

An \emph{affine ruling} (a $\PP^1$-ruling, a \emph{$\C^{(n*)}$-ruling}) is a morphism from a surface onto a smooth curve with general fiber isomorphic to $\C^1$ (respectively to $\PP^1$, $\C^1$ with $n$ points deleted). All these morphisms are called rational rulings.

\bdfn Let $p\:X\to B$ be a rational ruling  of a normal surface $X$. A triple $(\ov X,D,\ov p)$ is called \emph{a completion of $p$} if and only if $(\ov X,D)$ is a normal completion of $X$ and $\ov p:\ov X\to \ov B$ is a $\PP^1$-ruling onto a smooth curve $\ov B\supseteq B$ extending $p$. Such a completion is \emph{minimal} (we say also that $(\ov X,D)$ is $\ov p$-minimal) if it does not dominate any other completion of $p$. \edfn

If an snc-divisor $T$ (or rather its dual graph) is a chain and $T=T_1+T_2+\ldots+T_n$ is its decomposition into irreducible components so that $T_i\cdot T_{i+1}=1$ for $i=1,\ldots,n-1$, we then write $T=[-T_1^2,\ldots,-T_n^2]$. As long as $T$ is not considered as a twig attached to some other divisor containing $T$ there is no preferred choice of the tip ($T_1$ or $T_n$) of $T$, so in this case $T=[-T_n^2,\ldots,-T_1^2]$ as well. If $T$ is a twig of some fixed bigger divisor then by convention we always choose the tip of $T$ which is a tip of this bigger divisor as $T_1$.

The Kodaira dimension of a divisor $F$ on a smooth complete surface (hence projective by the theorem of Zariski) $\ov X$ is defined as $$\kappa(F)=\sup_{n>0} \dim \im(\Phi_{|nF|})\in \{-\8,0,1,2\},$$ where $\Phi_{|nF|}\:X\map \PP^N$ is the mapping given by the linear system $|nF|$. Then the logarithmic Kodaira dimension of a smooth open surface $X$ can be defined by taking some smooth completion $(\ov X,D)$ of $X$ and by putting $$\ovk(X)=\kappa(K_{\ov X}+D),$$ where $K_{\ov X}$ is a canonical divisor on $\ov X$. This is well-known to be independent of the smooth completion (see \cite{Iitaka} for the properties of $\kappa(F)$). The Kodaira dimension of a singular surface is defined to be the Kodaira dimension of any resolution. If $R$ is an snc-divisor on a complete surface and $Q(R)$ is the intersection matrix of $R$ then we define the \emph{discriminant} of $R$ by $d(R)=d(-Q(R))$.

\section{Balanced and standard completions}

One of the basic steps to take when dealing with an open surface is to construct a completion and boundary. Of course, these are not unique, but we want to bring to the attention of the reader some normalizing conditions, which make them more unique and more useful in practice. We came to this problem when trying to distinguish (or to find an isomorphism) between some $\Q$-acyclic surfaces. This type of analysis was done at least partially by many authors, with the most complete treatment in terms of weighted dual graphs in \cite{Daigle-weighted-graphs} and \cite{FKZ-weighted-graphs}. Here we present the necessary results using partially our own terminology. For simplicity we restrict ourselves to divisors whose dual graphs contain no loops (forests).

\bdfn\label{def:balanced} A rational chain $D=[a_1,\ldots,a_n]$ is \emph{balanced} if $a_1,\ldots,a_n\in\{0,2,3,\ldots\}$ or if $D=[1]$. A reduced snc-forest $D$ is \emph{balanced} if all rational chains contained in $D$ which do not contain branching components of $D$ are balanced. A normal pair $(X,D)$ is \emph{balanced} if $D$ is balanced. \edfn

The word \emph{balanced} stands here for the property that on one hand we do not allow non-branching $(-1)$-curves, but on the other hand we do not allow non-branching $b$-curves with positive $b$. The following operation is responsible for non-uniqueness of balanced completions of a given surface.

\bdfn\label{def2:flow} Let $(X,D)$ be a normal pair. Let $L$ be a $0$-curve which is a non-branching component of $D$. Make a blowup of a point $c\in L$, such that $c\in L\cap (D-L)$ if $L\cdot (D-L)=2$ and contract the proper transform of $L$. The resulting pair $(X',D')$, where $D'$ is the reduced direct image of the total transform of $D$ is called an \emph{elementary transform of $(X,D)$}. The point $c\in L$ is the \emph{center} of the transformation. A composition of elementary transformations of $D$ and its subsequent elementary transforms is called \emph{a flow inside $D$}. \edfn

For example, taking $T=[0,0,a_1,\ldots,a_n]$ one can obtain by a flow exactly the chains of type $[0,b,a_1,\ldots,a_n]$, $[a_1,\ldots,a_{k-1},a_k-b,0,b,a_{k+1},\ldots,a_n]$ or $[a_1,\ldots,a_n,b,0]$ where $1\leq k\leq n$ and $b\in\Z$. In particular, it is easy to see when two rational chains differ by a flow. The following result is the main property of balanced completions, it follows from \cite[3.36]{FKZ-weighted-graphs}.
\
\bprop\label{prop:boundary uniqueness} Any normal surface which admits a normal completion with a forest as a boundary has a balanced completion. Two such completions differ by a flow inside the boundary. In particular, all balanced boundaries of a given surface are isomorphic as curves. \eprop

Of course, two balanced boundaries of a given surface are in general non-isomorphic as weighted curves (weights are here the self-intersections of their components). One introduces normalizing conditions to deal with this. Following \cite{FKZ-weighted-graphs} the components of the divisor which remain after subtracting all non-rational and all branching components of an snc-divisor $D$ will be called the \emph{segments of $D$}. An snc-divisor is \emph{standard} if and only if any of its connected components is either $[1]$ or has all segments of type $[0]$, $[0,0]$, $[0,0,0]$, $[0,0,a_1,\ldots,a_n]$ or $[a_1,\ldots,a_n]$ with $a_1,\ldots,a_n\geq 2$. It follows from the example above and from the Hodge index theorem that each balanced chain can be carried to a standard form by a flow.

\bex The affine plane $\C^2$ has a balanced completion $(\F_2,T)$, where $\F_2$ is the second Hirzebruch surface and $T=[0,2]$. It follows that each balanced boundary of $\C^2$ is of type $[0,a]$ for some $a\geq 2$, it is standard if and only if $a=0$. Note that the isomorphism type of the boundary can change when we admit non-balanced completions, for example the boundary of $\C^2$ embedded in $\PP^2$ is $[-1]$. \eex

More generally, it follows from the above considerations about a flow inside chains that if a surface has a completion with some rational chain as a boundary then there are at most two weighted graphs which can be dual graphs of a standard boundary of this surface. If there are two then one of them is of type $T=[0,0,a_1,\ldots,a_n]$ with some $a_1,\ldots,a_n\geq 2$ and the second is \emph{the reversion} of the first one, $T^{rev}=[0,0,a_n,\ldots,a_1]$. Given one weighted dual graph of some standard boundary of a given surface each other can be easily described, which gives a refined method of distinguishing between many open surfaces.

\section{Basic properties}

Let $R$ be a ring. An \emph{$R$-homology plane} is a normal surface $X$ with $H^*(X,R)\cong R$. (This is a bit nonstandard, as usually $R$-homology planes are defined as smooth by definition). We say that $X$ is a \emph{homotopy plane} if $\pi_i(X)\cong 0$ for $i>0$. A $\Z$-homology plane with trivial $\pi_1$ is a homotopy plane by the theorem of Hurewicz. By a theorem of Whitehead homotopy planes are contractible.

\bex\label{ex:C2/G} Let $G<GL(2,\C)$ be a finite subgroup without pseudoreflections. Then $\C^2/G$ is a singular logarithmic homotopy plane for which the smooth locus has negative Kodaira dimension. Indeed, the linear contraction of $\C^2$ to $0\in \C^2$ descends to a contraction of the quotient and the quotient morphism $\C^2-\{0\}\to \C^2-\{0\}/G$ is \'{e}tale, so $\ovk(\C^2-\{0\}/G)=\ovk(\C^2-\{0\})=-\8$. \eex

We now fix the notation for the rest of the paper. Let $\epsilon\:S\to S'$ be a minimal snc-resolution of a singular $\Q$-homology plane $S'$. Denote the exceptional divisor of $\epsilon$ by $\E$. Let $(\ov S,D)$ be a smooth completion of $S$ and let $S_0$ be the smooth locus of $S'$. Since $S'$ is normal, its singular locus consists of a finite number of points. For topological spaces $A\subseteq X$ we write $H_i(X,A)$ for $H_i(X,A,\Q)$ and $b_i(X,A)$ for $\dim H_i(X,A)$.

Logarithmic $\Q$-homology planes are known to be affine by an argument of Fujita (cf. \cite[2.4]{Fujita}). They are also known to be rational due to Gurjar-Pradeep-Shastri \cite{PS-rationality}, \cite{GPS-logQhp_rational}, \cite{GP}. On the other hand, the following example shows that non-logarithmic $\Q$-homology planes can be non-rational.

\bex\label{ex:cone over a curve} Let $C\subseteq \PP^n$ be a projectively normal embedding of a smooth projective curve. Then the affine cone over $C$ is normal and contractible. It has a standard cylinder resolution, for which the exceptional divisor is isomorphic to $C$. In case $C$ is not rational, the cone is non-rational and non-logarithmic.\eex

In general we have the following result (cf. \cite[3.2, 3.4]{Palka-classification}).

\bthm Every singular $\Q$-homology plane is affine and birationally equivalent to a product $B\times \C^1$ for some curve $B$. \ethm

As for the affiness we note that the above mentioned argument of Fujita works if one assumes only the inclusion $\iota\:D\cup \E \to \ov S$ induces an isomorphism $H_2(\iota)\:H_2(D\cup \E)\to H_2(\ov S)$. One can prove that this condition is always satisfied. If $\E$ is a rational forest this can be seen as follows. First note that $H_i(\ov S,D\cup \E)\cong H^{4-i}(S_0)$ by the Lefschetz duality. One can prove that $b_2(S_0)=b_1(S_0)=b_1(\E)$, so if $\E$ is a rational forest then $b_3(\ov S,D\cup \E)=b_2(\ov S,D\cup \E)=0$ and $H_2(\iota)$ is an isomorphism. For $b_1(\E)\neq 0$ the argument is more complicated. For the rest we refer to loc. cit.

\bcor If $S'$ is a singular $\Q$-homology plane then its boundary is connected and the homology groups $H_i(S',\Z)$ vanish for $i\geq 2$. \ecor

\begin{proof} Since $S'$ is affine, by \cite{Karchjauskas_Letschetz_thm} its boundary is connected and $S'$ is homotopy equivalent to a CW-complex of real dimension at most two, hence $H_2(S',\Z)$ is torsionfree and $H_3(S',\Z)=H_4(S',\Z)=0$. Since $b_2(S')=0$, we get $H_2(S',\Z)=0$. \end{proof}

The rationality of smooth $\Q$-homology planes has strong consequences for vector bundles.

\bthm\label{thm:vb_over_Qhp} Let $S'$ be a smooth $\Q$-homology plane. Then any vector bundle over $S'$ is a sum of a line bundle and a trivial vector bundle. Moreover, $\Pic S'\cong H_1(X,\Z)$, so if $S'$ is a $\Z$-homology plane then all vector bundles over $S'$ are trivial. \ethm

\begin{proof} The first part of the theorem is a result of Murthy \cite{Murthy-vb_on_affine_ruled_surf} valid for smooth affine surfaces birational to a product of a curve with an line. Since $S'$ is rational, we have $\Pic S'\cong H^2(S',\Z)$. The groups $H_1(S',\Z)$ and $H_2(S',\Z)$ are finite, so by the universal coefficient theorem $H^2(S',\Z)\cong H_1(S',\Z)$. \end{proof}

We note that the question if a vector bundle over any smooth contractible threefold is necessarily trivial is open.

\section{Case $\ovk(S_0)=-\8$ and moduli}\label{sec:negative k(S_0)}

Let us start with recalling some known results in case $S_0$ has negative Kodaira dimension. This is a strong assumption, which in particular forces $S'$ to be logarithmic (and hence rational). This can be seen as follows.

Suppose $S_0$ is affine-ruled. The affine ruling extends to a $\PP^1$-ruling of some smooth completion $(\ov S,D+\E)$ of $S_0$, where $\E$ is the exceptional divisor of some (not necessarily minimal) snc-resolution of singularities of $S'$. The unique section of this extension contained in the boundary is in fact contained in $D$. Indeed, if it is contained in $\E$, then $D$ is vertical for this ruling (i.e. all its components intersect trivially with fibers), so since the homology classes of components of $D+\E$ generate $H_2(\ov S)$, the intersection form on $\ov S$ is semi-negative definite, contradicting the Hodge index theorem. Thus the affine ruling of $S_0$ extends to an affine ruling of $S'$, so $S'$ has at most cyclic singularities by \cite[Theorem 1]{Miyan-cylinderlike}. There is no bound on the number of singularities, $H_1(S',\Z)$ can be any finite abelian group. The extension has a unique fiber contained in the boundary (which can be assumed to be smooth) and each singular fiber contains a unique component not contained in $D\cup\E$ (see \cite[\S1, \S2]{MiSu-hPlanes} for more details). Note that from 1.2 loc. cit. if follows that $\C^2$ is the only smooth $\Z$-homology plane of negative Kodaira dimension (see \cite{Fuj-Zar} for the first proof of this) and this characterization implies the positive solution of the two-dimensional Cancellation Conjecture.

Now if $S_0$ is not affine-ruled (this can happen only if $S'$ is singular) then it follows from an important structure theorem by Miyanishi-Tsunoda \cite{MiTs-PlatFibr} that it contains an open subset $U$ with a very special $\C^*$-ruling called a \emph{Platonic fibration}. In fact one shows that $S_0=U$ (cf. \cite[3.1]{KR-ContrSurf}), which implies that $S'\cong \C^2/G$ for a finite small non-cyclic subgroup $G<GL(2,\C)$.

It is well known that there are only finitely many \emph{fake projective planes}, i.e. smooth projective surfaces with Betti numbers of $\PP^2$ (cf. \cite[V.1]{BHPV}). This is clearly not the case for $\Q$-homology planes (the name \emph{fake affine planes} is not used). Moreover, there are arbitrarily high dimensional nontrivial families of $\Q$-homology planes with a common weighted boundary. The following example of a family of singular $\Q$-homology planes is a modification of a similar example in \cite[4.16]{FZ-deformations}. We sketch the arguments.

\begin{figure}[h]\centering\includegraphics[scale=0.4]{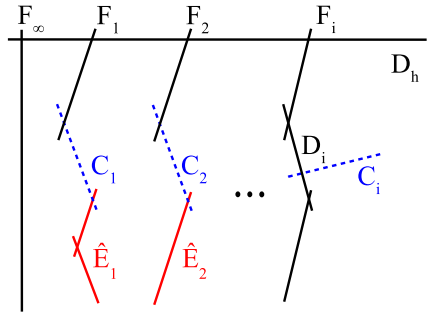}\caption{Singular fibers in the example \ref{ex:moduli}} \label{fig:moduli example}\end{figure}

\bex\label{ex:moduli} Let $p\:\F_1\to B\cong \PP^1$ be the $\PP^1$-ruling of the first Hirzebruch surface and let $N$ be a positive integer. Choose $N+3$ distinct points $x_\8,x_1,x_2,\ldots,x_{N+2}$ on the negative section. Blow up successively over each $x_i$ so to produce singular fibers with reductions $F_\8=[0]$, $F_1=[3,1,2,2]$, $F_2=[2,1,2]$, $F_i'=[2,1,2]$ for $i=2,\ldots,N+2$ lying respectively over $x_\8, x_1, \ldots, x_{N+2}$. Denote the resulting surface by $V$ and the proper transform of the negative section by $D_h$. Note there are chains $\E_1=[2,2]$ and $\E_2=[2]$ contained in $F_1$ and $F_2$ respectively, where we can assume that $\E_1, \E_2$ do not intersect $D_h$. For each $i\in \{3,\ldots,N\}$ choose a point $y_i$ on the $(-1)$-curve $D_i$ of $F_i'$ and blow up once. For $i\in\{1,\ldots,N\}$ denote the unique $(-1)$-curve of $F_i$ by $C_i$. Let $\ov S_{y}$, where $y=(y_3,\ldots,y_{N+2})$ be the resulting surface and put $D=F_\8+D_h+ (F_1-C_1-\E_1)+ (F_2-C_2-\E_2)+ \sum_{i=3}^N(F_i-C_i)$. One checks that the surface $S'_y$ obtained by the contraction of $\E_1$ and $\E_2$ on $\ov S-D$ is a singular $\Q$-homology plane. Clearly, the above family $S'_y$ is $N$-dimensional. Now if $S'_y\cong S'_z$ then the isomorphism lifts to $\ov S-D$ and then, since $F_\8$ is the only $0$-curve in $D$, extends to $\ov S_y-F_\8\cong \ov S_z-F_\8$ by \ref{prop:boundary uniqueness}, which in turn descends to an automorphism of $U=\F_1-F_\8-D_h\cong \C^2$ fixing fibers. However, if $x, y$ are respectively the horizontal and vertical coordinate on $U$, each automorphism of $U$ fixing fibers can be written as $(x,y)\to (x,\lambda y+P(x))$ for some $P[x]\in \C[x]$ and its lifting to $V$ acts by $\lambda^2$ on $D_i$ in some coordinates on $D_i$ (the multiplicity of $D_i$ in the fiber is $2$). Thus if we consider an $(N-1)$-dimensional subfamily with fixed $y_3\in D_3$ then $\lambda ^2=1$, so the mentioned action on each $D_i$ is trivial, hence $y=z$ and different members of this subfamily are non-isomorphic. \eex

\section{Logarithmic Bogomolov-Miyaoka-Yau}

In studies on singular $\Q$-homology planes we use often the logarithmic version of the Bogomolov-Miyaoka-Yau inequality proved by Kobayashi (cf. \cite{Kob}). Usually this inequality is stated for surfaces of general type in terms of the so-called \emph{strongly minimal model}. This is not necessary and in fact the inequality works for surfaces of non-negative Kodaira dimension. For example, the following lemma has been proved in \cite[2.5]{Palka-exceptional} as an easy corollary from an inequality of Bogomolov-Miyaoka-Yau type proved by Langer (cf. \cite{Langer}), which in particular generalizes both the inequality of Kobayashi and another inequality of Miyaoka \cite{Miyaoka}. For the notion of the Zariski decomposition and of the bark $\Bk D$ of an effective snc-divisor $D$ see \cite[\S 2.3]{Miyan-OpenSurf}. Put $D^\#=D-\Bk D$.

\blem\label{lem:KobIneq} Let $(X,D)$ be a smooth pair with $\kp(K_X+D)\geq 0$. Then:\benum[(i)]

\item $$3\chi(X-D)+\frac{1}{4}((K_X+D)^-)^2\geq (K_X+D)^2.$$

\item For each connected component of $D$, which is a connected component of $\Bk D$ (hence contractible to a quotient singularity) denote by $G_P$ the local fundamental group of the respective singular point $P$. Then $$\chi(X-D)+\sum_P\frac{1}{|G_P|}\geq \frac{1}{3} (K_X+D^\#)^2.$$ \eenum \elem

\noin We now illustrate the usefulness of the second inequality in our context.

\bcor\label{lem3:general is logarithmic}Let $S_0$ be the smooth locus of a singular $\Q$-homology plane $S'$ and let $\E_1, \ldots, \E_q$ be the connected components of the exceptional divisor of the snc-minimal resolution.  \benum[(i)]

\item If $\ovk(S_0)=2$ then $S'$ is logarithmic and $q=1$.

\item If $\ovk(S_0)=0$ or $1$ then either $q=1$, or $q=2$ and $\E_1=\E_2=[2]$.\eenum\ecor

\begin{proof}
Let $(S_m,D_m)$ be the almost minimal model of $(\ov S,D+\E)$. Since $S'$ is affine, $S_m-D_m$ is isomorphic to an open subset of $S_0$ satisfying $\chi(S_m-D_m)\leq \chi(S_0)=1-q$. Let $Q$ be the set of singular points which have been created by contracting connected components of $D$ as in \ref{lem:KobIneq}(ii). Since $\#Q\leq q$, the above inequality gives $$\frac{1}{3}((K_{S_m}+D_m)^+)^2\leq \chi(S_m-D_m)+\sum_{P\in Q}\frac{1}{|G_P|}\leq 1-q+\frac{\#Q}{2}\leq 1-\frac{q}{2}.$$

Now if $\ovk(S_0)=2$ then $((K_{S_m}+D_m)^+)^2>0$ and we get $q=1$ and $0<\sum_{P\in Q}\frac{1}{|G_P|}$, so there is a unique singular point on $S'$ and it is of quotient type. If $q>1$ and $\ovk(S_0)=0,1$ then $((K_{S_m}+D_m)^+)^2=0$ and we get $q=2$ and $1=1/|G_{P_1}|+1/|G_{P_2}|$, so $|G_{P_1}|=|G_{P_2}|=2$.
\end{proof}

\section{Exceptional $\Q$-homology planes}

The structure theorem for $\Q$-homology planes with smooth locus of non-general type is based on general structure theorems for open surfaces.

\bdfn A $\Q$-homology plane for which the smooth locus is neither of non-general type, nor $\C^1$- or $\C^*$-ruled is \emph{exceptional}. \edfn

Now the mentioned structure theorems lead to the fact that for exceptional $S'$ one has $\ovk(S_0)=0$. It was proved by Fujita (cf. \cite[8.64]{Fujita}) that each exceptional smooth $\Q$-homology plane is up to isomorphism one of three surfaces called $Y\{3,3,3\}$, $Y\{2,4,4\}$ and $Y\{2,3,6\}$ (the Fujita's surfaces of type $H[k,-k]$ with $k\geq 1$ are $\C^*$-ruled). The snc-minimal boundary of $Y\{a,b,c\}$ is a rational fork (a tree with one branching component) and its three maximal twigs have discriminants equal to $a,b,c$ respectively. A singular exceptional $\Q$-homology plane having a boundary with this property will be denoted by $SY\{a,b,c\}$. Together with the description of exceptional singular $\Q$-homology planes in \cite{Palka-exceptional} we have the following theorem (to be precise one still needs to prove that smooth $Y\{a,b,c\}$'s are not $\C^*$-ruled, but this can be done as in loc. cit.).

\bthm\label{thm:structure thm for non general type} If $S'$ is a $\Q$-homology plane with smooth locus $S_0$ of non-general type then $S_0$ is affine-ruled or $\C^*$-ruled or $S'$ is up to isomorphism one of five exceptional $\Q$-homology planes having smooth locus of Kodaira dimension zero: smooth $Y\{3,3,3\}$, $Y\{2,4,4\}$, $Y\{2,3,6\}$ and singular $SY\{3,3,3\}$, $SY\{2,4,4\}$. For the last two surfaces $\ovk(S')=0$ and the singular locus consists of a unique point of Dynkin type $A_2$ and $A_1$, respectively.\ethm

None of the exceptional surfaces is a $\Z$-homology plane. The maximal twigs of their snc-minimal boundaries consist of $(-2)$-curves. Writing this review the author noticed that a posteriori there is another (different than the one discovered in loc. cit.) nice description of exceptional singular $\Q$-homology planes. Namely, $Y\{2,4,4\}$ and $Y\{3,3,3\}$ have automorphism groups $\Z_2$ and $\Z_3$ respectively (the automorphism group of $Y\{2,3,6\}$ is trivial) and the actions have unique fixed points. The quotients are two non-isomorphic $\Q$-homology planes with smooth loci of Kodaira dimension zero. Suppose, say, the quotient $S'=Y\{3,3,3\}/\Aut Y\{3,3,3\}$ is not exceptional. Then its smooth locus $S_0$ is $\C^*$-ruled, so the Stein factorization of the pull-back of this $\C^*$-ruling gives a $\C^*$-ruling of the complement of the fixed point of $Y\{3,3,3\}$. Since $Y\{3,3,3\}$ is exceptional, the closures of the fibers meet in the fixed point, hence the closures of the fibers of the $\C^*$-ruling of $S_0$ meet in the singular point of $S'$. Thus the last $\C^*$-ruling does not extend to a ruling of $S'$. Since the singularity is cyclic, it follows from 5.4 \cite{Palka-classification} that $\ovk(S_0)=-\8$, a contradiction. (One can also get a contradiction with the fact that $Y\{3,3,3\}$ does not contain infinitely many contractible curves, cf. \ref{thm:contractible curves}).

We recall here the construction of $SY\{3,3,3\}$, mainly because of its beautiful connection with classical geometry. A \emph{projective configuration of type $(a_c,b_d)$} is an arrangement of $b$ lines in a projective space and $a$ points on these lines, such that each point belongs to $c$ lines and each line contains $d$ points. Clearly, $ac=bd$ for such a configuration.

\begin{figure}[h]\centering\includegraphics[scale=0.45]{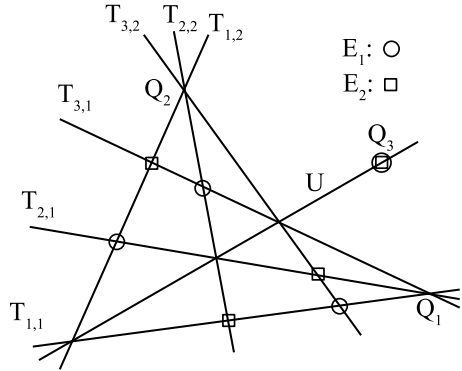}\caption{Singular Y\{3,3,3\}, dual Hesse configuration.}  \label{fig:Y333}\end{figure}

\bex\label{ex:Y333} Up to a projective automorphism there exists a unique projective configuration $\mathfrak{H}$ of type $(12_3,9_4)$ (the uniqueness is easy to show using information on the automorphism group of the configuration which we give below). The dual configuration $\frak H^*$ is the famous Hesse configuration $(9_4,12_3)$ of flexes of an elliptic curve. It is known that $Aut(\frak H)\cong Aut(\frak H^*)$ has order $216$ and is isomorphic to the group of special affine transformations of $\F_3^2$, i.e. $\F_3^2\rtimes SL(2,\F_3)$ (cf. \cite[\S4]{Dolgachev-abstract_configurations}). Choose three points $Q_1,Q_2,Q_3$ of $\mathfrak{H}$ none two of which lie on a common line of $\frak H$ and choose some line $U$ incident to $Q_3$. By taking the dual choice in $\frak H$ we check using linear algebra that $Aut(\frak H)$ acts transitively on the set of such choices and the stabilizer $\Gamma<Aut(\frak H)$ has order three. For $i=1,2$ let $T_{1,i}$, $T_{2,i}$, $T_{3,i}$ be the lines incident to $Q_i$ and let $E_1,E_2,U$ be the lines incident to $Q_3$. Blow up once in each of eight points not incident to $U$. Let $\ov S$ be the resulting complete surface. We denote divisors and their proper transforms by the same letters. Let $B$ be the exceptional curve over $Q_2$, put $D=B+T_{1,1}+T_{1,2}+T_{2,1}+T_{2,2}+T_{3,1}+T_{3,2}$ and $\E=E_1+E_2$. Clearly, all components of $D-B+\E$ are $(-2)$-curves, $D$ and $\E$ are disjoint and $D$ is a fork. Put $S'=(\ov S-D)/\E$. One can show that $Aut(S')\cong \Gamma$.

To see that $S'$ is $\Q$-acyclic note that $b_1(\ov S)=0$ and $b_2(\ov S)=b_2(D\cup\E)=9$. Now since $d(D+\E)\neq 0$, the natural morphism $H_2(D\cup \E)\to H_2(\ov S)$ is an isomorphism. The homology exact sequence of the pair $(\ov S,D)$ and Lefschetz duality give $b_1(S)=b_3(S)=b_4(S)=0$ and $b_2(S)=\#\E$ ($S=\ov S\setminus D$). We know from the above that $H_2(\E)\to H_2(S)$ is a monomorphism, so the homology exact sequence of the pair $(S,\E)$ gives that $S'$ is $\Q$-acyclic.

We check easily that in both cases $K_{\ov S} + D^\#=K_{\ov S}+D^\#+\E^\#$ intersects trivially with all components of $D+\E$, hence $\ovk(S_0)=\ovk(S')=0$. See \cite[\S 5]{Palka-exceptional} for an explicit realization of $\frak H$ and for a proof that the constructed $\Q$-homology plane does not admit a $\C^*$-ruling.\eex

\section{Non-logarithmic $\Q$-homology planes}

Let us start with a generalization of example \ref{ex:cone over a curve}.

\bex\label{ex:cone/G} Let $U$ be an affine cone over a projectively normal curve and let $G$ be a finite group acting on it so that the vertex of the cone is the unique fixed point, the action is free on its complement and recpects lines of the cone. Then the quotient $S'=U/G$ is a normal contractible surface. \eex

What is surprising is that all non-logarithmic $\Q$-homology planes arise in the above way. Namely, we have the following theorem (see \cite[5.8]{Palka-classification}).

\bthm\label{thm:nonlogarithmic are cones/G} Every singular $\Q$-homology plane containing a non-quotient singularity is a quotient of an affine cone over a smooth projective curve by an action of a finite group which is free off the vertex of the cone and respects the set of lines through the vertex. In particular, it is contractible, has negative Kodaira dimension, has a unique singular point and its smooth locus is $\C^*$-ruled. \ethm

\begin{figure}[h]\centering\includegraphics[scale=0.4]{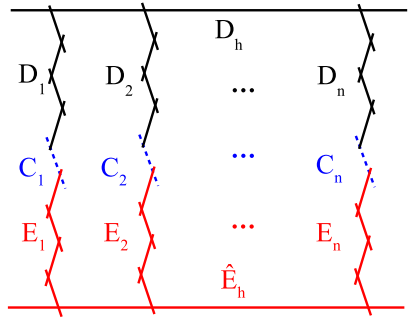}\caption{Construction of non-logarithmic $\Q$-homology planes.}  \label{fig:nonextendable}\end{figure}

\begin{proof}(sketch) By \ref{lem3:general is logarithmic} and \ref{thm:structure thm for non general type} the smooth locus of a non-logarithmic $\Q$-homology plane $S'$ is $\C^1$- or $\C^*$-ruled. By the results of section \ref{sec:negative k(S_0)}, $\ovk(S_0)\geq 0$, so in fact $S_0$ is $\C^*$-ruled. One can show (cf. 3.6 loc. cit.) that if this ruling extends to a $\C^*$-ruling of $S'$ then $S'$ is necessarily logarithmic, so we can further assume that this is not the case. This means that there is a $\PP^1$-ruling $p:\ov S\to \PP^1$ of some completion $(\ov S,D+\E)$ of $S_0$ as before, such that for a general fiber $f$ we have $f\cdot D=f\cdot \E=1$. We can assume that the completion is $p$-minimal. Then each singular fiber of $p$ is a so-called \emph{columnar fiber}, which means that it is a rational snc-chain of discriminant zero, its components have negative self-intersections, it contains a unique $(-1)$-curve $C_i$ which is also a unique $S_0$-component of the fiber (i.e. a component not contained in $D\cup \E$) and it is intersected by the horizontal components of $D$ and $\E$ in tips. It follows that $D$ and $\E$ contain unique branching components $D_h$ and $\E_h$. Write the singular fibers as $D_i+C_i+E_i$, where $D_i\subseteq D$ and $E_i\subseteq \E$  (see Fig.~\ref{fig:nonextendable}). Contracting the singular fibers we get a $\PP^1$-bundle over a complete curve, so all non-logarithmic singular $\Q$-homology planes can be reconstructed starting from such a bundle by producing columnar singular fibers, taking out $D$ and contracting $\E$. The mentioned bundle admits a usual $\C^*$-action fixing pointwise the images of $D_h$ and $\E_h$. This action induces a $\C^*$-action on $S'$ with the singular point as the unique fixed point. Thus by \cite[1.1]{Pinkham-C*-surfaces} $S'$ is a quotient of an affine cone over a smooth projective curve by an action of a finite group. In fact without knowing this global description the contractibility follows also from \cite[5.9, 4.19]{Fujita} and the Whitehead theorem, as one can show that $\pi_1(\E_h)\to\pi_1(S)$ is an isomorphism and $\pi_1(S)\to \pi_1(S')$ is an epimorphism. Since $S=\ov S-D$ is affine-ruled, we get $\ovk(S')=-\8$. \end{proof}

We see that one can obtain any non-logarithmic $\Q$-homology plane by creating its snc-minimal completion $\ov S$ of the resolution as in the proof above by starting with a non-trivial $\PP^1$-bundle over some smooth complete curve $B$ and by blowing up some number $n$ of columnar fibers. Then exactly one of two forks, say $\E$, separated by vertical $(-1)$-curves is negative definite, hence by \cite{Grauert} contracts analytically to a singularity. What is quite surprising here is that this contraction is always algebraic (cf. the contraction criterion \cite[2.6]{Palka-classification}). Suppose $\E_h\cong B\cong \PP^1$. Although $\E$ is a rational tree, the singularity does not have to be a rational singularity. Indeed, using Artin's criterion one can show that if $\E_h^2+n\leq 0$ then $S'$ has a rational and if $\E_h^2+n\geq 2$ a non-rational singularity (cf. \cite[5.8]{Palka-classification}). The case $\E_h^2+n=1$ is more subtle (cf. \cite[5.8]{Pinkham-C*-surfaces}).

\section{$\C^*$-ruled $\Q$-homology planes}

By the results described in previous sections the classification of $\Q$-homo\-logy planes of non-general type reduces now to the classification of logarithmic (and hence rational) $\Q$-homology planes which are $\C^*$-ruled, or in other words, for which the smooth locus is $\C^*$-ruled and the ruling extends to a $\C^*$-ruling of $S'$. Note that the existence of a $\C^*$-ruling implies that $\ovk(S_0)\neq 2$ by the 'easy addition theorem' (cf. \cite[2.1.5]{Miyan-OpenSurf}). This case was analyzed in \cite{MiSu-hPlanes}, where one can find a description of singular fibers and the computation of $\ovk(S_0)$ and $H_1(S',\Z)$ in terms of these fibers. We note here that in \cite[6.8]{Palka-classification} we have redone some incorrect computations of $\ovk(S_0)$ from loc. cit. (identified then with the Kodaira dimension of $S'$) and we have computed $\ovk(S')$ too. We discuss two issues here.

First, it is practically useful to know when a $\C^{(n*)}$-ruled surface is a $\Q$-homology plane. To formulate a criterion we need to recall the definition of some numbers characterizing rational rulings. Having a fixed $\PP^1$-ruling of a smooth complete surface $X$ and a reduced divisor $T$ we define  $$\Sigma_{X-T}=\underset{F\nsubseteq T}{\sum}(\sigma(F)-1),$$ where $\sigma(F)$ is the number of $(X-T)$-components (i.e. irreducible components not contained in $T$) of a fiber $F$ (cf. \cite[4.16]{Fujita}). $T_h$ is the horizontal part of $T$, which consists of components of $T$ intersecting nontrivially with a general fiber. If $T_h=0$ then $T$ is \emph{vertical}. The numbers $h$ and $\nu$ are defined respectively as $\#T_h$ and as the number of fibers contained in $T$.

\blem\label{lem:when ruling gives sQhp} Let $(\ov S,T)$ be a smooth pair and let $p\:\ov S\to \PP^1$ be a $\PP^1$-ruling. Assume the following conditions are satisfied:\benum[(i)]

\item there exists a unique connected component $D$ of $T$ which is not vertical,

\item $D$ is a rational tree,

\item $\Sigma_{\ov S-T}=h+\nu-2$,

\item $d(D)\neq 0$.

\eenum Then the surface $S'$ defined as the image of $\ov S-D$ after contraction of connected components of $T-D$ to points is a rational $\Q$-homology plane and $p$ induces a rational ruling of $S'$. Conversely, if $p'\:S'\to B$ is a rational ruling of a rational $\Q$-homology plane $S'$ (not necessarily singular) then any completion $(\ov S,T,p)$ of the restriction of $p'$ to the smooth locus of $S'$ has the above properties. \elem

The conditions (iii)-(iv) are equivalent to the fact that $H_2(D\cup\E)\to H_2(\ov S)$ is an isomorphism, similar criteria were used by many authors. What is important, in case of a $\C^*$-ruling the most problematic condition (iv) can be replaced by an easier and more geometric condition (cf. \cite[6.1]{Palka-classification}).

Second, there is a question about the uniqueness of a $\C^*$-ruling of $S'$. Let us assume $S'$ is singular. In case $\ovk(S_0)=1$ it is easy to prove that there is a unique $\C^*$-ruling of $S'$ and it is induced by the $\C^*$-ruling of $S_0$ given by the linear system of some multiple of the logarithmic canonical divisor of $S_0$. In case $\ovk(S_0)=0$ generically there are two $\C^*$-rulings of $S'$, but there may be zero (exceptional $\Q$-homology planes), one or three as well (cf. 6.12 loc. cit.). We will see that this information is important for example when computing the number of contractible curves on $S'$.

As for the completions of $S'$, a non-affine-ruled $S'$ has a unique balanced completion unless it admits an untwisted $\C^*$-ruling with base $\C^1$ (untwisted means that for some completion of the surface the ruling extends to a $\PP^1$-ruling for which the horizontal part of the boundary consists of two irreducible components). In the latter case there are exactly two balanced completions (cf. 6.11 loc. cit.).

\section{Contractible curves}\label{sec:contractible curves}

It is well known that the logarithmic Bogomolov-Miyaoka-Yau inequality imposes restrictions on the number of contractible curves on $S'$. According to the author's knowledge this number is known except the cases when $S'$ is non-logarithmic or when $S'$ is singular and $\ovk(S_0)=0$ (see \cite{Zaid_isotrivial}, \cite{Zaid_Corr_isotrivial} \cite{GM-Affine-lines}, \cite{GP-affine_lines_on_Qhp}). In the non-logarithmic case it is easily seen to be infinity by \ref{thm:nonlogarithmic are cones/G} and, as we sketch below, in the last case it can be deduced from the knowledge on the number (and types) of $\C^*$-rulings of $S_0$ \cite[6.12]{Palka-classification}. The final result is as follows.

\bthm\label{thm:contractible curves} Let $l$ be the number of contractible curves on a $\Q$-homology plane $S'$. Let $S_0$ be the smooth locus of $S'$. Then: \benum[(i)]

\item if $\ovk(S_0)=2$ then $l=0$,

\item if $S'$ is exceptional (hence $\ovk(S')=\ovk(S_0)=0$) then $l=0$,

\item if $S'$ is non-logarithmic (hence $\ovk(S')=-\8$, $\ovk(S_0)=0,1$) then $l=\8$,

\item if $\ovk(S_0)=-\8$ then $l=\8$,

\item in other cases ($S'$ is $\C^*$-ruled and $\ovk(S_0)=0,1$), $l=1$ or $2$.

\eenum \ethm

\begin{proof}(sketch) If $\ovk(S_0)=-\8$ then $S'$ is affine-ruled or isomorphic to $\C^2/G$ for some $G<GL(2,\C)$, so $l=\8$. If $S'$ is non-logarithmic then by \ref{thm:nonlogarithmic are cones/G} $l=\8$. We can therefore assume that $\ovk(S_0)\geq 0$ and $S'$ is logarithmic. Suppose $S'$ contains a contractible curve $L$. It follows from the logarithmic Bogomolov-Miyaoka-Yau inequality that $\ovk(S_0-L)\leq 1$ (cf. \cite{GM-Affine-lines}). Since $S'$ is rational, $Pic(S_0)=\Coker (\Pic (D+\E)\to \Pic \ov S)$ is torsion, so the class of $L$ in $\Pic(S_0)$ is torsion. Then there is a morphism $f\:S_0-L\to \C^*$ and taking its Stein factorization one gets a $\C^*$-ruling of $S_0-L$, which (as $\ovk(S_0)\geq 0$) extends to a $\C^*$-ruling of $S_0$. Since $S'$ is logarithmic, each $\C^*$-ruling of $S_0$ extends in turn to a $\C^*$-ruling of $S'$. Therefore, any contractible curve on $S'$ is vertical for some $\C^*$-ruling of $S'$. In particular, if $l>0$ then $S'$, and hence $S_0$, is necessarily $C^*$-ruled, hence $S'$ cannot be exceptional. The analysis of fibers (Suzuki's formula) leads to the corollary that there exist one or two contractible vertical curves for a given $\C^*$-ruling. Now if $\ovk(S_0)=1$ then this ruling is unique (given by a multiple of the logarithmic canonical divisor of $S_0$), hence $l=1,2$. Consider now the case $\ovk(S_0)=0$. Here the problem is more difficult, as there may be more $\C^*$-rulings of $S'$. In case $S'$ is smooth it was shown in \cite{GP-affine_lines_on_Qhp} that $\l=1$. In case $S'$ is singular we have computed the number and types of possible $\C^*$-rulings of $S_0$ in \cite[6.12]{Palka-classification} (one can do the same if $S'$ is smooth in a similar way). Since this number is finite we see that $l$ is finite and nonzero. Looking more closely at the proof one can show that $l\leq 2$. \end{proof}

\bex Let $\pi\:S\to \C^2$ be the restriction of the projection $(x,y,z)\to (x,y)$ to the surface $S=\{(x,y,z)\in \C^3: z^n=f(x,y)\}$, where $n$ is a positive natural number. Then $\pi$ is a branched cover with the curve $C=\{(x,y)\in \C^2:f(x,y)=0\}$ as the branch locus. Suppose $S$ is a smooth $\Q$-homology plane. Then $C$ is smooth and we have $1-\chi(C)=1-\chi(\pi^{-1}(C))=1-\chi(S)+\chi(S-\pi^{-1}(C))=n\chi(\C^2-C)=n(1-\chi(C))$, hence $\chi(C)=1$. If $U$ is any smooth affine curve then it is non-complete, so its Euler characteristic is smaller than the Euler characteristic of the smooth completion $\overline{U}$, hence $\chi(U)\leq 1-2g(\overline{U})$, where $g$ is the genus. In case $U$ is irreducible it follows that $\chi(U)\geq 0$ only if $U\cong \C^1$ or $U\cong \C^*$. In particular, some component of $C\cong \pi^{-1}(C)$ is an affine line, so by \ref{thm:contractible curves} $S$ is not of general type. In fact either $\ovk(S)=-\8$ or $S$ is $\C^*$-ruled. Smooth $\Q$-homology planes of this kind have been classified in \cite{Maharana-Qhp_cyclic_covers}.  \eex

\section{Smooth $\Q$-homology planes of general type}\label{sec:general type}

Let $S'$ be a smooth $\Q$-homology plane of general type, i.e. $\ovk(S')=2$. We know already that $S'$ is rational. By \ref{thm:contractible curves} $S'$ contains no contractible curves, which implies that its snc-minimal completion $(\ov S,D)$ is almost minimal. We can assume that this completion is balanced. Since $S'$ is neither $\C^1$- nor $\C^*$-ruled, $D$ contains no non-branching $0$-curves, so any flow inside $D$ is trivial. It follows that the balanced completion of $S'$ is unique. As for now there is no classification, but there are some partial results.

The first example of such a surface which was shown to be non-isomorphic to $\C^2$ is the famous example of Ramanujam (\cite{Ramanujam}). We know now infinitely many examples of this kind. One method of construction of $\Q$-homology planes of general type is to use $\C^{(n*)}$-rulings with $n\geq 2$. For $n=2$ this was done in \cite{MiSu-Qhp's_withCstarstar}. Another method was used by tom Dieck and Petrie.

\bdfn Let $(\ov S,D)$ be a completion of $S'$ for which there exists a birational morphism $f:(\ov S,D)\to (\PP^2,f_*D)$. Then $f_*D$ is called a \emph{plane divisor} of $S'$. If $f_*D$ is a sum of lines then it is a \emph{linear plane divisor} and we say that $S'$ comes from the \emph{line arrangement} $f_*D$. \edfn

In \cite{Dieck_linear_div_hp} it was noticed that an inequality of Hirzebruch bounds the number of types of possible linear plane divisors for smooth $\Q$-homology planes, a list of these divisors has been given. In \cite{DieckPetrie-hp_and_alge_curves} a general algorithm for recovering smooth $\Q$-homology planes (in fact countable series of them) starting from a given rational divisor on a minimal rational complete surface is described. The conjecture that all smooth $\Q$-homology planes of general type have linear plane divisors is not true by an example of tom Dieck (cf. \cite{Dieck_symmetric}). Tom Dieck's smooth $\Z$-homology plane has a nontrivial automorphism group (it is necessarily finite, as $S'$ is of general type), so is also a counterexample to the earlier conjecture of Petrie \cite{Petrie_Aut_of_affine_surfaces}.

\bex Let $C\subseteq \PP^2$ be an irreducible curve. Recall that a singular point $p\in C$ is a \emph{cusp} if $C$ is locally irreducible at $p$. Since $C\hookrightarrow \PP^2$ induces a monomorphism on $H_2(\cdot ,\Q)$, by the Lefschetz duality and by the long exact sequence of the pair $(\PP^2,C)$ the Betti numbers of $S_C=\PP^2\setminus C$ are $b_2(S_C)=b_1(C)$ and $b_i(S_C)=0$ for $i=1$ and $i>2$. Assume that $C$ is rational and cuspidal. Then $S_C$ is a smooth $\Q$-homology plane (this is in fact equivalent). By \cite{Wakab-kod_of_curve_complement} if $C$ has two cusps then $\ovk(S_C)\geq 0$ and if it has more than two cusps then $S_C$ is of general type. The literature on plane cuspidal curves is rich, see for example \cite{FZ-cuspidal}, \cite{BLMN-rat_cusp_plane_curves} and references there. \eex

We now list some conjectural properties of smooth $\Q$-homology planes.

\begin{conjecture} Let $S'$ be a smooth $\Q$-homology plane of general type and let $(\ov S,D)$ be its minimal smooth completion. \benum[(A)]

\item $S'$ has a plane divisor consisting of lines and conics.

\item $S'$ admits a $\C^{(3*)}$-ruling.

\item $(K_{\ov S}+D)^2=-2$, or equivalently $K_{\ov S}\cdot(K_{\ov S}+D)=0$.

\item $S'$ is rigid and has unobstructed deformations.

\item The set of all possible Eisenbud-Neumann diagrams for smooth $\Q$-homology planes is finite.

\eenum\end{conjecture}

In \cite{Dieck-Optimal_homotopy} conjectures (A)-(C) have been stated and verified for all known smooth $\Z$-homology planes. Sugie \cite{Sugie_C3-4star_fibred_qhp's} analyzed $\C^{(n*)}$-rulings on smooth $\Q$-homology planes and classified possible singular fibers for $n=2$. Flenner and Zaidenberg \cite[6.12]{FZ-deformations} have shown that part (D), which implies (C), holds for smooth $\Q$-homology planes of general type having linear plane divisors. See \cite{Zaidenberg-open_problems} for (D), (E) and related conjectures.

Recently, the following result was proved \cite{GKMR}.

\bthm Let $S'$ be a smooth homotopy plane of general type. Then $\Aut S'$ is cyclic, its action on $S'$ has a unique fixed point and is free off this point.\ethm

The situation for $\Aut S'$ of a $\Q$-homology plane is not known.

\section{Smooth locus of general type}

Let $S'$ be a singular $\Q$-homology plane with smooth locus $S_0$ of general type, i.e. $\ovk(S_0)=2$. By \ref{lem3:general is logarithmic} $S'$ has a unique singular point and this point is of quotient type. This means that $\E$ is either a chain if this point is a cyclic singularity or a negative definite rational fork if not. Since $S'$ is logarithmic, it is rational. By \ref{thm:contractible curves} it contains no contractible curves, which implies that the snc-minimal completion $(\ov S,D+\E)$ of $S_0$ is almost minimal. Arguing as above we see that in fact this snc-minimal completion is unique and does not contain non-branching $0$-curves. Again, there is no classification, but there are some partial results.

Note that due to the existence of the transfer homomorphism for branched coverings (cf. \cite{Smith_transfer}) the quotient of a smooth $\Q$-homology plane of general type by its automorphism group is a $\Q$-homology plane (with smooth locus of general type).

A priori there is no restriction on the Kodaira dimension of $S'$. However, refining the methods of Koras-Russell \cite{KR-ContrSurf} M. Koras and the author have obtained the following theorem (this is also a part of the thesis of the author written under the supervision of M. Koras, cf. \cite{PaKo-general_type}).

\bthm Singular $\Q$-homology planes with smooth locus of general type have non-negative Kodaira dimension. \ethm

The cases $\ovk(S')=0,1$ and especially the case $\ovk(S')=2$ have not been analyzed systematically as for now. Also the number of known examples is much smaller than in the case of smooth $\Q$-homology planes. For example, similarly as in the smooth case, one could ask if $S_0$ has a completion $(\ov S,D+\E)$ admitting a birational morphism $f\:\ov S\to \PP^2$ with $f_*(D+\E)$ a sum of lines and conics. Even the case when $f_*(D+\E)$ is a sum of lines has not been studied.

The known examples of singular $\Z$-homology planes of general type the local fundamental group of the singular point (which has order equal to $d(\E)$) is cyclic of order not bigger than six. As for now the following result in this direction has been proven (\cite{GKMR}):

\bthm A singular $\Z$-homology plane with smooth locus of general type has a cyclic quotient singularity. \ethm

\bibliographystyle{amsalpha}
\bibliography{bibl}

\providecommand{\bysame}{\leavevmode\hbox to3em{\hrulefill}\thinspace}
\providecommand{\MR}{\relax\ifhmode\unskip\space\fi MR }
\providecommand{\MRhref}[2]{%
  \href{http://www.ams.org/mathscinet-getitem?mr=#1}{#2}
}
\providecommand{\href}[2]{#2}
\begin{thebibliography}{FdBLMHN07}

\bibitem[BHPVdV04]{BHPV}
Wolf~P. Barth, Klaus Hulek, Chris A.~M. Peters, and Antonius Van~de Ven,
  \emph{Compact complex surfaces}, second ed., Ergebnisse der Mathematik und
  ihrer Grenzgebiete. 3. Folge. A Series of Modern Surveys in Mathematics
  [Results in Mathematics and Related Areas. 3rd Series. A Series of Modern
  Surveys in Mathematics], vol.~4, Springer-Verlag, Berlin, 2004.

\bibitem[Dai03]{Daigle-weighted-graphs}
Daniel Daigle, \emph{Classification of weighted graphs up to blowing-up and
  blowing-down}, \href{http://arxiv.org/abs/math/0305029}{arXiv:math/0305029}
  (2003).

\bibitem[Dol04]{Dolgachev-abstract_configurations}
Igor~V. Dolgachev, \emph{Abstract configurations in algebraic geometry}, The
  {F}ano {C}onference, Univ. Torino, Turin, 2004, pp.~423--462,
  (\href{http://arxiv.org/abs/math/0304258}{arXiv:math/0304258}).

\bibitem[FdBLMHN07]{BLMN-rat_cusp_plane_curves}
J.~Fern{\'a}ndez~de Bobadilla, I.~Luengo, A.~Melle-Hern{\'a}ndez, and
  A.~N{\'e}methi, \emph{On rational cuspidal plane curves, open surfaces and
  local singularities}, Singularity theory, World Sci. Publ., Hackensack, NJ,
  \href{http://arxiv.org/abs/math/0604421}{arXiv:math/0604421}, 2007,
  pp.~411--442.

\bibitem[FKZ07]{FKZ-weighted-graphs}
Hubert Flenner, Shulim Kaliman, and Mikhail Zaidenberg, \emph{Birational
  transformations of weighted graphs}, Affine algebraic geometry, Osaka Univ.
  Press, Osaka, 2007, pp.~107--147.

\bibitem[Fuj79]{Fuj-Zar}
Takao Fujita, \emph{On {Z}ariski problem}, Proc. Japan Acad. Ser. A Math. Sci.
  \textbf{55} (1979), no.~3, 106--110.

\bibitem[Fuj82]{Fujita}
\bysame, \emph{On the topology of noncomplete algebraic surfaces}, J. Fac. Sci.
  Univ. Tokyo Sect. IA Math. \textbf{29} (1982), no.~3, 503--566.

\bibitem[FZ94]{FZ-deformations}
Hubert Flenner and Mikhail Zaidenberg, \emph{{$\bold Q$}-acyclic surfaces and
  their deformations}, Classification of algebraic varieties ({L}'{A}quila,
  1992), Contemp. Math., vol. 162, Amer. Math. Soc., Providence, RI, 1994,
  pp.~143--208.

\bibitem[FZ96]{FZ-cuspidal}
\bysame, \emph{On a class of rational cuspidal plane curves}, Manuscripta Math.
  \textbf{89} (1996), no.~4, 439--459.

\bibitem[GKMR10]{GKMR}
R.~V. Gurjar, M.~Koras, M.~Miyanishi, and P.~Russell, \emph{A homology plane of
  general type can have at most a cyclic quotient singularity}, preprint, 62
  pages, 2010.

\bibitem[GM92]{GM-Affine-lines}
R.~V. Gurjar and M.~Miyanishi, \emph{Affine lines on logarithmic {${\bf
  Q}$}-homology planes}, Math. Ann. \textbf{294} (1992), no.~3, 463--482.

\bibitem[GP95]{GP-affine_lines_on_Qhp}
R.~V. Gurjar and A.~J. Parameswaran, \emph{Affine lines on {${\bf Q}$}-homology
  planes}, J. Math. Kyoto Univ. \textbf{35} (1995), no.~1, 63--77.

\bibitem[GP99]{GP}
R.~V. Gurjar and C.~R. Pradeep, \emph{{${\bf Q}$}-homology planes are rational.
  {III}}, Osaka J. Math. \textbf{36} (1999), no.~2, 259--335.

\bibitem[GPS97]{GPS-logQhp_rational}
R.~V. Gurjar, C.~R. Pradeep, and Anant.~R. Shastri, \emph{On rationality of
  logarithmic {$\bold Q$}-homology planes. {II}}, Osaka J. Math. \textbf{34}
  (1997), no.~3, 725--743.

\bibitem[Gra62]{Grauert}
Hans Grauert, \emph{\"{U}ber {M}odifikationen und exzeptionelle analytische
  {M}engen}, Math. Ann. \textbf{146} (1962), 331--368.

\bibitem[Iit82]{Iitaka}
Shigeru Iitaka, \emph{Algebraic geometry}, Graduate Texts in Mathematics,
  vol.~76, Springer-Verlag, New York, 1982, An introduction to birational
  geometry of algebraic varieties, North-Holland Mathematical Library, 24.

\bibitem[Kar77]{Karchjauskas_Letschetz_thm}
K.~K. Kar{\v{c}}jauskas, \emph{A generalized {L}efschetz theorem}, Funkcional.
  Anal. i Prilo\v zen. \textbf{11} (1977), no.~4, 80--81.

\bibitem[Kob90]{Kob}
Ryoichi Kobayashi, \emph{Uniformization of complex surfaces}, K\"ahler metric
  and moduli spaces, Adv. Stud. Pure Math., vol.~18, Academic Press, Boston,
  MA, 1990, pp.~313--394.

\bibitem[KR07]{KR-ContrSurf}
Mariusz Koras and Peter Russell, \emph{Contractible affine surfaces with
  quotient singularities}, Transform. Groups \textbf{12} (2007), no.~2,
  293--340.

\bibitem[Lan03]{Langer}
Adrian Langer, \emph{Logarithmic orbifold {E}uler numbers of surfaces with
  applications}, Proc. London Math. Soc. (3) \textbf{86} (2003), no.~2,
  358--396.

\bibitem[Mah09]{Maharana-Qhp_cyclic_covers}
Alok Maharana, \emph{{$Q$}-homology planes as cyclic covers of {$\Bbb A\sp
  2$}}, J. Math. Soc. Japan \textbf{61} (2009), no.~2, 393--425.

\bibitem[Miy81]{Miyan-cylinderlike}
Masayoshi Miyanishi, \emph{Singularities of normal affine surfaces containing
  cylinderlike open sets}, J. Algebra \textbf{68} (1981), no.~2, 268--275.

\bibitem[Miy84]{Miyaoka}
Yoichi Miyaoka, \emph{The maximal number of quotient singularities on surfaces
  with given numerical invariants}, Math. Ann. \textbf{268} (1984), no.~2,
  159--171.

\bibitem[Miy01]{Miyan-OpenSurf}
Masayoshi Miyanishi, \emph{Open algebraic surfaces}, CRM Monograph Series,
  vol.~12, American Mathematical Society, Providence, RI, 2001.

\bibitem[Miy07]{Miyan_recent_developments}
\bysame, \emph{Recent developments in affine algebraic geometry: from the
  personal viewpoints of the author}, Affine algebraic geometry, Osaka Univ.
  Press, Osaka, 2007, pp.~307--378.

\bibitem[MS91a]{MiSu-hPlanes}
M.~Miyanishi and T.~Sugie, \emph{Homology planes with quotient singularities},
  J. Math. Kyoto Univ. \textbf{31} (1991), no.~3, 755--788.

\bibitem[MS91b]{MiSu-Qhp's_withCstarstar}
Masayoshi Miyanishi and Tohru Sugie, \emph{{${\bf Q}$}-homology planes with
  {${\bf C}\sp {**}$}-fibrations}, Osaka J. Math. \textbf{28} (1991), no.~1,
  1--26.

\bibitem[MT84]{MiTs-PlatFibr}
Masayoshi Miyanishi and Shuichiro Tsunoda, \emph{Noncomplete algebraic surfaces
  with logarithmic {K}odaira dimension {$-\infty$} and with nonconnected
  boundaries at infinity}, Japan. J. Math. (N.S.) \textbf{10} (1984), no.~2,
  195--242.

\bibitem[Mur69]{Murthy-vb_on_affine_ruled_surf}
M.~Pavaman Murthy, \emph{Vector bundles over affine surfaces birationally
  equivalent to a ruled surface}, Ann. of Math. (2) \textbf{89} (1969),
  242--253.

\bibitem[Pal09]{Palka-classification}
Karol Palka, \emph{On the classification of singular {$\bold Q$}-acyclic
  surfaces}, \href{http://arxiv.org/abs/0806.3110}{arXiv:0806.3110} (2009).

\bibitem[Pal10]{Palka-exceptional}
\bysame, \emph{Exceptional singular {$\bold Q$}-homology planes},
  \href{http://arxiv.org/abs/0909.0772}{arXiv:0909.0772}, to appear in Annales
  Inst. Fourier (2010).

\bibitem[Pet89]{Petrie_Aut_of_affine_surfaces}
Ted Petrie, \emph{Algebraic automorphisms of smooth affine surfaces}, Invent.
  Math. \textbf{95} (1989), no.~2, 355--378.

\bibitem[Pin77]{Pinkham-C*-surfaces}
H.~Pinkham, \emph{Normal surface singularities with {$C\sp*$} action}, Math.
  Ann. \textbf{227} (1977), no.~2, 183--193.

\bibitem[PK10]{PaKo-general_type}
Karol Palka and Mariusz Koras, \emph{Singular {$\bold Q$}-homology planes of
  negative kodaira dimension have smooth locus of non-general type},
  \href{http://arxiv.org/abs/1001.2256}{arXiv:1001.2256} (2010).

\bibitem[PS97]{PS-rationality}
C.~R. Pradeep and Anant~R. Shastri, \emph{On rationality of logarithmic {$\bold
  Q$}-homology planes. {I}}, Osaka J. Math. \textbf{34} (1997), no.~2,
  429--456.

\bibitem[Ram71]{Ramanujam}
C.~P. Ramanujam, \emph{A topological characterisation of the affine plane as an
  algebraic variety}, Ann. of Math. (2) \textbf{94} (1971), 69--88.

\bibitem[Smi83]{Smith_transfer}
Larry Smith, \emph{Transfer and ramified coverings}, Math. Proc. Cambridge
  Philos. Soc. \textbf{93} (1983), no.~3, 485--493.

\bibitem[Sug99]{Sugie_C3-4star_fibred_qhp's}
Toru Sugie, \emph{Singular fibers of homology planes with pencils of rational
  curves}, Mem. Fac. Educ. Shiga Univ. III Nat. Sci. (1999), no.~49, 29--40
  (2000).

\bibitem[tD90a]{Dieck_linear_div_hp}
Tammo tom Dieck, \emph{Linear plane divisors of homology planes}, J. Fac. Sci.
  Univ. Tokyo Sect. IA Math. \textbf{37} (1990), no.~1, 33--69.

\bibitem[tD90b]{Dieck_symmetric}
\bysame, \emph{Symmetric homology planes}, Math. Ann. \textbf{286} (1990),
  no.~1-3, 143--152.

\bibitem[tD92]{Dieck-Optimal_homotopy}
\bysame, \emph{Optimal rational curves and homotopy planes}, Bol. Soc. Mat.
  Mexicana (2) \textbf{37} (1992), no.~1-2, 115--138, Papers in honor of
  Jos{\'e} Adem (Spanish).

\bibitem[tDP93]{DieckPetrie-hp_and_alge_curves}
Tammo tom Dieck and Ted Petrie, \emph{Homology planes and algebraic curves},
  Osaka J. Math. \textbf{30} (1993), no.~4, 855--886.

\bibitem[Wak78]{Wakab-kod_of_curve_complement}
Isao Wakabayashi, \emph{On the logarithmic {K}odaira dimension of the
  complement of a curve in {$P^{2}$}}, Proc. Japan Acad. Ser. A Math. Sci.
  \textbf{54} (1978), no.~6, 157--162.

\bibitem[Za{\u\i}87]{Zaid_isotrivial}
M.~G. Za{\u\i}denberg, \emph{Isotrivial families of curves on affine surfaces,
  and the characterization of the affine plane}, Izv. Akad. Nauk SSSR Ser. Mat.
  \textbf{51} (1987), no.~3, 534--567, 688.

\bibitem[Za{\u\i}91]{Zaid_Corr_isotrivial}
\bysame, \emph{Additions and corrections to the paper: ``{I}sotrivial families
  of curves on affine surfaces, and the characterization of the affine plane''
  [{I}zv.\ {A}kad.\ {N}auk {SSSR} {S}er.\ {M}at.\ {\bf 51} (1987), no.\ 3,
  534--567]}, Izv. Akad. Nauk SSSR Ser. Mat. \textbf{55} (1991), no.~2,
  444--446.

\bibitem[Za{\u\i}99]{Zaid_exotic_lectures}
M.~Za{\u\i}denberg, \emph{Exotic algebraic structures on affine spaces},
  Algebra i Analiz \textbf{11} (1999), no.~5, 3--73,
  \href{http://arxiv.org/abs/math/9801075}{arXiv:9801075}.

\bibitem[Za{\u\i}05]{Zaidenberg-open_problems}
Mikhail Za{\u\i}denberg, \emph{Selected problems},
  \href{http://arxiv.org/abs/math/0501457}{arXiv:0501457} (2005).

\end{thebibliography}
\end{document}